\documentclass[11pt]{article}
\usepackage[T1]{fontenc}
\usepackage[utf8]{inputenc}
\usepackage{amsmath,amsthm,amscd,amsfonts,amssymb,graphicx,color}
\usepackage[margin=1in]{geometry}
\usepackage[hidelinks]{hyperref}

\newtheorem{thm}{Theorem}[section]

\numberwithin{equation}{section}
\newcommand{\wt}{\mathrm{wt}}
\newcommand{\tvs}{\mathrm{tvs}}
\newcommand{\bv}[1]{\mathbf{#1}}

\title{A Unified Approach to Total Vertex Irregularity Strength}
\author{Aleams Barra\\
Algebra Research Group, Institut Teknologi Bandung, Indonesia\\
\texttt{aleamsbarra@itb.ac.id}}
\date{}

\begin{document}
\maketitle

\begin{abstract}
We develop a unified framework for computing the total vertex irregularity strength (tvs) of diverse graph classes, significantly extending our prior work in Barra et al. This comprehensive approach generalizes and enhances earlier methods, offering a robust technique applicable to a wide range of graph structures. A notable achievement is resolving an open problem by determining the tvs for simple 2-regular graphs, confirming the conjecture by Ahmad et al. Furthermore, our framework provides unified proofs for the tvs of cycles, paths, prisms, wheels, helm graphs, and friendship graphs. Importantly, the versatility of this method suggests its potential utility in computing tvs for additional graph classes, broadening its impact in graph theory and computational applications.
\end{abstract}

\section{Introduction}
A network can be modeled as a graph, where multiple connections between vertices are represented by multiple edges. To represent the number of edges between vertices, Chartrand et al.~\cite{chartrand1988irregular} proposed replacing multiple edges with a single edge labeled with a positive integer.

The \emph{vertex sum} of a vertex in such a graph is defined as the sum of the labels on all edges incident to that vertex. Chartrand et al.~further introduced the concept of an \emph{irregular network}, characterized by a graph in which all vertices have distinct vertex sums. 

This framework defines the \emph{irregularity strength} of a graph \(G\), denoted \(s(G)\), as the smallest positive integer \(k\) for which there exists an edge labeling of \(E(G)\) with \(\{1, 2, \ldots, k\}\) ensuring distinct vertex sums.

For \(d\)-regular graphs of order \(n\), one trivially has \(s(G)\ge n/d\). Faudree and Lehel conjectured the existence of an absolute constant \(c\) such that \(s(G)\le n/d + c\). A strong asymptotic form proved by Przyby{\l}o and Wei~\cite{przybylo2023asymptotic} states that for any fixed \(\varepsilon\in(0,1/4)\), there exist constants \(c_1,c_2\) such that for graphs on \(n\) vertices with minimum degree \(\delta>1\) and no isolated edge,
\[
s(G)\le \frac{n}{\delta}\left(1+\frac{c_1}{\delta^{\varepsilon}}\right)+c_2.
\]
This result refines earlier bounds, including the general estimate of Kalkowski et al.~\cite{kalkowski2011upper} and dense-graph bounds of Majerski and Przyby{\l}o~\cite{majerski2014dense}; a short follow-up~\cite{przybylo2023short} further provides explicit additive constants in regular-graph regimes.

Baca et al.~\cite{bacaIrregularTotalLabellings2007} extended the concept of irregularity strength \(s(G)\) to its \emph{total-labeling version}, where both vertices and edges of \(G\) are assigned positive integers.
The \emph{vertex sum} of a vertex \(v \in V(G)\) is given by
\[
\omega_E(v) := \sum_{vw \in E(G)} \lambda(vw),
\]
where the sum includes all edges incident to \(v\). The \emph{weight} of \(v\) is then defined as
\[
\mathrm{wt}(v) = \omega_E(v) + \lambda(v).
\]

A total labeling \(\lambda\) is deemed \(k\)-irregular if \(\mathrm{range}(\lambda) \subseteq \{1, 2, \ldots, k\}\) and all vertices in \(G\) have distinct weights. The smallest positive integer \(k\) for which \(\lambda\) is \(k\)-irregular is termed the \emph{total vertex irregularity strength} of \(G\), denoted \(\mathrm{tvs}(G)\).

Exact values of \(\mathrm{tvs}(G)\) have been established for several graph classes, including cycles, paths, prisms, trees, wheels, helm graphs, friendship graphs, complete graphs, and complete bipartite graphs (see~\cite{bacaIrregularTotalLabellings2007, nurdinTotalVertexIrregularity2010, wijaya2008total, ahmadTotalVertexIrregularity, Wijaya2005}).

For general graphs, a widely used upper bound due to Anholcer et al.~\cite{anholcer2009upper} states that for every graph \(G\) of order \(n\) and minimum degree \(\delta>0\),
\[
\mathrm{tvs}(G)\le 3\left\lceil\frac{n}{\delta}\right\rceil+1.
\]
Further asymptotic improvements in dense-graph settings were later obtained by Majerski and Przyby{\l}o~\cite{majerski2013dense}.

Despite these developments, the standard approach to determining \(\mathrm{tvs}(G)\) still consists of two separate steps: establishing a lower bound and then proving its attainability via an explicit construction. While the lower-bound arguments are relatively consistent, the constructive part is often graph-specific, and no unified construction strategy has proved effective across diverse families. In many cases, the labeling itself is intricate and requires careful, lengthy verification.

Motivated by this gap, we introduce a construction technique for total vertex irregular labelings based on a shared underlying principle that applies to a broad range of graph classes.

\section{Preliminaries}
The concept of \(\{1,s\}\)-edge labeling was introduced by Barra et al.~\cite{barraTotalVertexIrregularity2025a} to determine the total vertex irregularity strength of cubic graphs with a perfect matching. Earlier work, including Baca et al.~\cite{bacaIrregularTotalLabellings2007}, also used two-value edge-label patterns (notably \(\{1,2\}\)) in specific constructions; our approach treats this idea systematically through a variable set \(\{1,s\}\). A \(\{1,s\}\)-edge labeling on a graph \(G\) is defined as a function
\[
\lambda: E(G) \to \{1, s\}.
\]
Using this labeling, the vertex sums of all vertices in \(G\) can be computed. This labeling can be extended to an \(s\)-irregular labeling by suitably defining
\[
\lambda: V(G) \to \{1, 2, \ldots, s\}
\]
such that, for any distinct vertices \(v, w \in V(G)\), the weights satisfy \(\mathrm{wt}(v) \neq \mathrm{wt}(w)\).

For a given \(\{1,s\}\)-edge labeling \(\lambda\) and a subset \(U \subseteq V(G)\), we denote by \(\mathbf{U}\) the vector of length \(|U|\), listing the vertices in \(U\) from left to right in non-decreasing order of their vertex sums.

For any positive integer \(a\), the boldface notation \(\mathbf{a}\) represents a constant vector with all entries equal to \(a\), where the length of the vector is implied by the context. Likewise, for positive integers \(a\) and \(b\) with \(a \leq b\), the notation \(\mathbf{[a,b]}\) denotes the vector
\[
\mathbf{[a,b]} = \langle a, a+1, \ldots, b \rangle.
\]

Given a vector \(\mathbf{U} = \langle u_1, u_2, \ldots, u_m \rangle\), we define
\begin{align*}
\omega_E(\mathbf{U}) &= \langle \omega_E(u_1), \omega_E(u_2), \ldots, \omega_E(u_m) \rangle, \\
\lambda(\mathbf{U}) &= \langle \lambda(u_1), \lambda(u_2), \ldots, \lambda(u_m) \rangle, \quad \text{and} \\
\mathrm{wt}(\mathbf{U}) &= \langle \mathrm{wt}(u_1), \mathrm{wt}(u_2), \ldots, \mathrm{wt}(u_m) \rangle.
\end{align*}

We frequently rely on the following lower bound, established in~\cite{bacaIrregularTotalLabellings2007}.

\begin{thm}
Let \(G\) be a graph with \(n\) vertices, where \(\delta\) and \(\Delta\) denote the minimum and maximum degrees of its vertices, respectively. Then
\begin{equation}\label{lbound}
\mathrm{tvs}(G) \geq \left\lceil \frac{n + \delta}{\Delta + 1} \right\rceil.
\end{equation}
\end{thm}
\section{Main Results}
In this section, we utilize the \(\{1,s\}\)-edge labeling method to compute several established results regarding the total vertex irregularity strength of various graphs. Our focus is on establishing a unified approach to determine the \(\mathrm{tvs}\) values of these graphs that simplifies their construction and enhances the transparency of verifying their correctness. While verifying the correctness of constructions in earlier studies can be intricate, our method ensures that the validation process is both straightforward and transparent.

As a culmination of this approach, we apply it at the end of this section to settle an open problem concerning the total vertex irregularity strength of simple 2-regular graphs.

\subsection{Cycles and Paths}

The construction for cycles presented here will later be adapted to determine the total vertex irregularity strength of simple 2-regular graphs.

\begin{thm}\label{cycle}
Let \(C_n\) be the cycle graph with \(n\) vertices. Then
\[
\mathrm{tvs}(C_n) = \lceil (n+2)/3 \rceil.
\]
\end{thm}

\begin{proof}
In \(C_n\), both the minimum degree \(\delta\) and maximum degree \(\Delta\) are 2. Applying the lower bound from~\eqref{lbound}, we obtain
\[
\mathrm{tvs}(C_n) \geq \left\lceil \frac{n + \delta}{\Delta + 1} \right\rceil = \left\lceil \frac{n + 2}{3} \right\rceil.
\]

Let \(n = a + 2b + c\), where \(a,c \geq 1\) and \(b \geq 0\). From an arbitrary starting edge, assign 1 to the first \(a\) consecutive edges of the cycle, label the next \(2b\) edges alternately with \(s\) and 1, and mark the remaining \(c\) edges with \(s\). 

This labeling yields distinct vertex sums: inner vertices of the 1-labeled segment have a vertex sum of 2 (being incident to two edges labeled 1), inner vertices of the alternating \((s,1)\)-labeled segment have a vertex sum of \(s + 1\), and inner vertices of the \(s\)-labeled segment have a vertex sum of \(2s\). Vertices at boundaries between segments also have a vertex sum of \(s + 1\).

Define \(V_2\), \(V_{s+1}\), and \(V_{2s}\) as the sets of vertices with vertex sums 2, \(s + 1\), and \(2s\), respectively. The distribution of vertex sums under \(\lambda\) is given by
\[
D_{\lambda}(G) = \left(|V_2|, |V_{s+1}|, |V_{2s}|\right).
\]
From our construction, this distribution is
\begin{equation}\label{distribution}
D_{\lambda} = (a - 1, 2b + 2, c - 1).
\end{equation}

Set \(s = \lceil (n + 2) / 3 \rceil\). We proceed with two cases based on the parity of \(s\):

\textbf{Case 1: \(s\) is even.} Choose \(a\) and \(b\) such that \(a - 1 = s - 1\) (i.e., \(a = s\)) and \(2b + 2 = s\) (i.e., \(b = (s - 2)/2\)). Then \(c = n - a - 2b = n + 2 - 2s\), yielding
\[
D_{\lambda} = (s - 1, s, n + 1 - 2s).
\]
Extend \(\lambda\) to \(V(C_n)\) by defining
\begin{equation}\label{s_even}
\lambda(\mathbf{V_2}) = \mathbf{[1,s-1]}, \quad \lambda(\mathbf{V_{s+1}}) = \mathbf{[1,s]}, \quad \text{and} \quad \lambda(\mathbf{V_{2s}}) = \mathbf{[2, n + 2 - 2s]}.
\end{equation}
The weights are then
\begin{align*}
\mathrm{wt}(\mathbf{V_2}) &= \mathbf{2} + \mathbf{[1,s-1]} = \mathbf{[3,s+1]}, \\
\mathrm{wt}(\mathbf{V_{s+1}}) &= \mathbf{s + 1} + \mathbf{[1,s]} = \mathbf{[s+2,2s+1]}, \\
\mathrm{wt}(\mathbf{V_{2s}}) &= \mathbf{2s} + \mathbf{[2,n+2-2s]} = \mathbf{[2s+2,n+2]}.
\end{align*}
These intervals are clearly pairwise disjoint. Moreover, since \(s=\left\lceil\frac{n+2}{3}\right\rceil\), we have \(n+2-2s\le s\). Hence all vertex labels are at most \(s\), and the resulting weights are distinct.

\textbf{Case 2: \(s\) is odd.} Choose \(a\) and \(b\) such that \(a - 1 = s\) (i.e., \(a = s + 1\)) and \(2b + 2 = s - 1\) (i.e., \(b = (s - 3)/2\)), so \(c = n + 2 - 2s\). This gives
\[
D_{\lambda} = (s, s - 1, n + 1 - 2s).
\]
Define
\begin{equation}\label{s_odd}
\lambda(\mathbf{V_2}) = \mathbf{[1,s]}, \quad \lambda(\mathbf{V_{s+1}}) = \mathbf{[2,s]}, \quad \text{and} \quad \lambda(\mathbf{V_{2s}}) = \mathbf{[2,n+2-2s]}.
\end{equation}
The weights become
\begin{align*}
\mathrm{wt}(\mathbf{V_2}) &= \mathbf{2} + \mathbf{[1,s]} = \mathbf{[3,s+2]}, \\
\mathrm{wt}(\mathbf{V_{s+1}}) &= \mathbf{s + 1} + \mathbf{[2,s]} = \mathbf{[s+3,2s+1]}, \\
\mathrm{wt}(\mathbf{V_{2s}}) &= \mathbf{2s} + \mathbf{[2,n+2-2s]} = \mathbf{[2s+2,n+2]}.
\end{align*}
Since \(n + 2 - 2s \leq s\) holds as before, all vertex labels are at most \(s\). The weight intervals are disjoint, ensuring all vertices have distinct weights.

In both cases, \(\lambda\) is \(s\)-irregular with \(s = \lceil (n + 2) / 3 \rceil\), matching the lower bound. Thus, \(\mathrm{tvs}(C_n) = \lceil (n + 2) / 3 \rceil\).
\end{proof}
Having established the result for cycles, we now apply a similar approach to compute the total vertex irregularity strength of paths
\begin{thm}
Let \(P_n\) be the path with \(n\geq 3\) vertices. Then
\[
\mathrm{tvs}(P_n) = \left\lceil \frac{n+1}{3} \right\rceil.
\]
\end{thm}

\begin{proof}
For \(P_n\), the minimum degree \(\delta = 1\) (at one endpoint) and maximum degree \(\Delta = 2\) (at inner vertices). Applying the lower bound from~\eqref{lbound}, we get
\[
\mathrm{tvs}(P_n) \geq \left\lceil \frac{n + \delta}{\Delta + 1} \right\rceil = \left\lceil \frac{n + 1}{3} \right\rceil.
\]

To prove this bound is achievable, we construct a total labeling that ensures all vertex weights are distinct while matching the bound. Since \(P_n\) has \(n - 1\) edges, let \(n - 1 = a + 2b + c\). Starting from the leftmost edge, define a \(\{1,s\}\)-edge labeling \(\lambda\): assign 1 to the first \(a\) consecutive edges, label the next \(2b\) edges alternately with \(s\) and 1, and label the remaining \(c\) edges with \(s\). We require \(a, c \geq 1\), while \(b \geq 0\).

The labeling yields vertex sums similar to the cycle case, with sums 1 at one endpoint, \(s\) at the other endpoint, 2 in the 1-labeled segment, \(s + 1\) in the alternating and boundary vertices, and \(2s\) in the \(s\)-labeled segment. Define \(V_2\) as the set of vertices with vertex sum \(1\) or \(2\), \(V_{s+1}\) as those with sum \(s\) or \(s + 1\), and \(V_{2s}\) as those with sum \(2s\). The distribution is
\[
D_{\lambda} = (|V_2|, |V_{s+1}|, |V_{2s}|) = (a, 2b + 1, c - 1).
\]

Set \(s = \left\lceil (n + 1) / 3 \right\rceil\). We consider two cases based on the parity of \(s\):

\textbf{Case 1: \(s\) is even.} Choose \(a\) and \(b\) such that \(a = s - 1\) and \(2b + 2 = s\). Then consequently, \(c = n - 1 - a - 2b = n + 2 - 2s\), so
\[
D_{\lambda} = (s - 1, s, n + 1 - 2s).
\]
Extend \(\lambda\) to \(V(P_n)\) as follows:
\[
\lambda(\bv{V_2}) = \langle 1, \bv{[1,s-2]} \rangle, \quad \lambda(\bv{V_{s+1}}) = \langle 1, \bv{[1,s-1]} \rangle, \quad \lambda(\bv{V_{2s}}) = \bv{[1, n+1-2s]},
\]
yielding weights
\begin{align*}
\wt(\bv{V_2}) &= \langle 1, \bv{2} \rangle + \langle 1, \bv{[1,s-2]} \rangle = \bv{[2,s]}, \\
\wt(\bv{V_{s+1}}) &= \langle s, \bv{s+1} \rangle + \langle 1, \bv{[1,s-1]} \rangle = \bv{[s+1,2s]}, \\
\wt(\bv{V_{2s}}) &= \bv{2s} + \bv{[1,n+1-2s]} = \bv{[2s+1,n+1]}.
\end{align*}
Since \(s = \lceil (n + 1) / 3 \rceil\), we have \(n + 1 \leq 3s\), so \(n + 1 - 2s \leq s\), ensuring all labels are at most \(s\). The weight intervals \(\bv{[2,s]}\), \(\bv{[s+1,2s]}\), and \(\bv{[2s+1,n+1]}\) are disjoint and cover all \(n\) vertices distinctly.

\textbf{Case 2: \(s\) is odd.} Choose \(a\) and \(b\) such that \(a = s\) and \(2b + 2 = s - 1\). Then consequently, \(c  = n + 2 - 2s\), so
\[
D_{\lambda} = (s, s - 1, n + 1 - 2s).
\]
Define
\[
\lambda(\bv{V_2}) = \langle 1, \bv{[1,s-1]} \rangle, \quad \lambda(\bv{V_{s+1}}) = \langle 2, \bv{[2,s-1]} \rangle, \quad \lambda(\bv{V_{2s}}) = \bv{[1,n+1-2s]},
\]
yielding
\begin{align*}
\wt(\bv{V_2}) &= \langle 1, \bv{2} \rangle + \langle 1, \bv{[1,s-1]} \rangle = \bv{[2,s+1]}, \\
\wt(\bv{V_{s+1}}) &= \langle s, \bv{s+1} \rangle + \langle 2, \bv{[2,s-1]} \rangle = \bv{[s+2,2s]}, \\
\wt(\bv{V_{2s}}) &= \bv{2s} + \bv{[1,n+1-2s]} = \bv{[2s+1,n+1]}.
\end{align*}
Again, \(n + 1 - 2s \leq s\) holds, and the weight intervals \(\bv{[2,s+1]}\), \(\bv{[s+2,2s]}\), and \(\bv{[2s+1,n+1]}\) are disjoint and cover all \(n\) vertices distinctly.

In both cases, \(\lambda\) is \(s\)-irregular with \(s = \lceil (n + 1) / 3 \rceil\), matching the lower bound. Thus, \(\mathrm{tvs}(P_n) = \lceil (n + 1) / 3 \rceil\).
\end{proof}

\subsection{Prisms and Wheels}
 The \emph{prism graph} \(D_n\) consists of two \(n\)-cycles, \(v_1 v_2 \cdots v_n v_1\) and \(w_1 w_2 \cdots w_n w_1\), together with the edges \(\{v_i w_i : i = 1, \ldots, n\}\) connecting corresponding vertices of the two cycles. The \emph{wheel graph} \(W_n\) consists of an \(n\)-cycle \(v_1 v_2 \cdots v_n v_1\) together with a center vertex \(x\), connected to each vertex of the cycle by the edges \(\{x v_i : i = 1, \ldots, n\}\). An edge of the form \(x v_i\) in \(W_n\) is called a \emph{spoke} of the wheel.

The total vertex irregularity strengths of prisms and wheels were first computed by Baca et al.~\cite{bacaIrregularTotalLabellings2007} and Wijaya et al.~\cite{wijaya2008total}, respectively. For wheel graphs, our construction simplifies the original approach and enhances the transparency of its verification.

\begin{thm}\label{prism}
Let \(D_n\) be the prism graph with \(n\) vertices in each cycle. Then
\[
\mathrm{tvs}(D_n) = \left\lceil \frac{n}{2} \right\rceil + 1.
\]
\end{thm}

\begin{proof}
The prism \(D_n\) has \(2n\) vertices, all of degree \(\delta = \Delta = 3\). Applying the lower bound from~\eqref{lbound}, we obtain
\[
\mathrm{tvs}(D_n) \geq \left\lceil \frac{2n + \delta}{\Delta + 1} \right\rceil = \left\lceil \frac{2n + 3}{4} \right\rceil = \left\lceil \frac{n}{2} \right\rceil + 1.
\]

To achieve this bound, we employ our \(\{1,s\}\)-labeling on \(D_n\). Set \(s = \left\lceil \frac{n}{2} \right\rceil + 1\). Label all edges of the cycle \(v_1 v_2 \cdots v_n v_1\) with 1 and all edges of the cycle \(w_1 w_2 \cdots w_n w_1\) with \(s\). For the cross edges \(v_i w_i\), label \(s - 1\) of them with 1 and the remaining \(n + 1 - s\) with \(s\).

This labeling produces vertex sums of 3 or \(s + 2\) for vertices \(v_i\) and \(2s + 1\) or \(3s\) for vertices \(w_i\), based on whether the cross edge \(v_i w_i\) is labeled 1 or \(s\). Define \(V_3\), \(V_{s+2}\), \(V_{2s+1}\), and \(V_{3s}\) as the sets of vertices with sums 3, \(s + 2\), \(2s + 1\), and \(3s\), respectively. The distribution is
\[
D_{\lambda} = (|V_3|, |V_{s+2}|, |V_{2s+1}|, |V_{3s}|) = (s - 1, n + 1 - s, s - 1, n + 1 - s).
\]

Extend \(\lambda\) to \(V(D_n)\) by defining
\[
\lambda(\bv{V_3}) = \lambda(\bv{V_{2s+1}}) = \bv{[1,s-1]}, \quad \lambda(\bv{V_{s+2}}) = \lambda(\bv{V_{3s}}) = \bv{[1,n+1-s]},
\]
yielding weights
\begin{align*}
\wt(\bv{V_3}) &= \bv{3} + \bv{[1,s-1]} = \bv{[4,s+2]}, \\
\wt(\bv{V_{s+2}}) &= \bv{s+2} + \bv{[1,n+1-s]} = \bv{[s+3,n+3]}, \\
\wt(\bv{V_{2s+1}}) &= \bv{2s+1} + \bv{[1,s-1]} = \bv{[2s+2,3s]}, \\
\wt(\bv{V_{3s}}) &= \bv{3s} + \bv{[1,n+1-s]} = \bv{[3s+1,n+2s+1]}.
\end{align*}
Since \(s = \lceil n/2 \rceil + 1 \geq (n/2) + 1\), we have \(n + 1 - s \leq s - 1\) and \(n + 3 \leq 2s + 1\), ensuring all labels are at most \(s\) and the intervals \(\bv{[4,s+2]}\), \(\bv{[s+3,n+3]}\), \(\bv{[2s+2,3s]}\), and \(\bv{[3s+1,n+2s+1]}\) are disjoint, covering \(2n\) distinct weights.

Thus, \(\lambda\) is \(s\)-irregular with \(s = \lceil n/2 \rceil + 1\), matching the bound, so \(\mathrm{tvs}(D_n) = \lceil n/2 \rceil + 1\).
\end{proof}

\begin{thm}\label{wheel}
Let \(W_n\) be the wheel graph with \(n + 1\) vertices, consisting of an \(n\)-cycle and a center vertex, for \(n \geq 3\). Then
\[
\mathrm{tvs}(W_n) = \lceil (n + 3) / 4 \rceil.
\]
\end{thm}

\begin{proof}
First, we establish a lower bound. For a cycle vertex \(v\) (degree 3), the smallest possible weight is 4 when all incident edges and \(v\) itself are labeled 1, while the largest possible weight is \(4k\) if \(k=\mathrm{tvs}(W_n)\), giving \(4k-3\) possible values for cycle-vertex weights. Since there are \(n\) cycle vertices, we need
\(
4k - 3 \geq n,
\)
so
\[
k \geq \lceil (n + 3) / 4 \rceil.
\]
Thus, \(\mathrm{tvs}(W_n) \geq \lceil (n + 3) / 4 \rceil\).

To attain this bound, we perform a \(\{1,s\}\)-edge labeling on \(W_n\). We label all edges in \(W_n\) (cycle edges and spokes) with 1 or \(s\). Excluding spoke labels, cycle vertices have partial vertex sums 2, \(s + 1\), or \(2s\); denote the corresponding sets by \(U_2\), \(U_{s+1}\), and \(U_{2s}\). Applying the cycle construction together with the general distribution formula~\eqref{distribution}, and choosing \(a=n+3-3s\), \(b=s-2\), and \(c=s+1\), we obtain
\[
|U_2| = n + 2 - 3s, \quad |U_{s+1}| = 2s - 2, \quad |U_{2s}| = s.
\]

Extend \(\lambda\) by labeling the spokes as follows: all edges from vertices in \(U_2\) to the center are labeled 1, all edges from vertices in \(U_{2s}\) to the center are labeled \(s\), and among the vertices in \(U_{s+1}\), exactly \(s-1\) spokes are labeled 1 while the remaining \(s-1\) spokes are labeled \(s\).

Then each vertex in \(U_2\) has full vertex sum \(2+1=3\). Since \(|U_2|=n+2-3s\), this gives \(|V_3|=n+2-3s\). Similarly, the \(s-1\) vertices of \(U_{s+1}\) whose spokes are labeled 1 have full vertex sum \(s+2\), so \(|V_{s+2}|=s-1\), while the remaining \(s-1\) vertices of \(U_{s+1}\) have full vertex sum \(2s+1\), so \(|V_{2s+1}|=s-1\). Finally, each vertex in \(U_{2s}\) has full vertex sum \(2s+s=3s\), and therefore \(|V_{3s}|=|U_{2s}|=s\). Hence the distribution of full vertex sums is
\[
(|V_3|, |V_{s+2}|, |V_{2s+1}|, |V_{3s}|) = (n + 2 - 3s, s - 1, s - 1, s).
\]

Label the cycle vertices:
\begin{gather*}
\lambda(\bv{V_3}) = \bv{[1,n+2-3s]}, \quad \lambda(\bv{V_{s+2}}) = \bv{[1,s-1]}, \\ \lambda(\bv{V_{3s}}) = \bv{[1,s]}, \quad \lambda(\bv{V_{2s+1}}) = \bv{[1,s-1]}.
\end{gather*}
The weights are
\begin{align*}
\wt(\bv{V_3}) &= \bv{3} + \bv{[1,n+2-3s]} = \bv{[4,n+5-3s]}, \\
\wt(\bv{V_{s+2}}) &= \bv{s+2} + \bv{[1,s-1]} = \bv{[s+3,2s+1]}, \\
\wt(\bv{V_{2s+1}}) &= \bv{2s+1} + \bv{[1,s-1]} = \bv{[2s+2,3s]}, \\
\wt(\bv{V_{3s}}) &= \bv{3s} + \bv{[1,s]} = \bv{[3s+1,4s]}.
\end{align*}
Set \(s = \lceil (n + 3) / 4 \rceil\). Then
\[
4s \geq n + 3 > n + 2 \Rightarrow n + 5 - 3s < s + 3.
\]
It follows that the weights of all cycle vertices are distinct.

For the center vertex \(x\), the number of spokes labeled \(s\) is \(2s-1\), while the remaining \(n-2s+1\) spokes are labeled 1. If we set \(\lambda(x)=s\), then
\[
\wt(x)=s(2s-1)+1(n-2s+1)+s=2s^2-2s+n+1.
\]
For \(n\ge 4\), we have \(s\ge 2\), and
\[
\wt(x)-4s=2s^2-6s+n+1 \ge 2s^2-6s+5 > 2s^2-6s+4 = 2(s-1)(s-2)\ge 0.
\]
Hence \(\wt(x)>4s\), so \(\wt(x)\) is distinct from all cycle-vertex weights. For \(n=3\), \(W_3\cong K_4\), and this case was already settled by Baca et al.~\cite{bacaIrregularTotalLabellings2007}.

\end{proof}
\subsection{Helm and Friendship Graphs}

 The \emph{helm graph} \(H_n\) is obtained from the \emph{wheel graph} \(W_n\) by attaching a pendant edge to each vertex of the \(n\)-cycle in \(W_n\), forming \(n\) pendant vertices and resulting in \(2n + 1\) vertices (\(n\) pendant vertices plus the \(n + 1\) vertices of \(W_n\)) and \(3n\) edges. The \emph{friendship graph} \(F_n\) is constructed by taking \(n\) copies of the triangle \(K_3\) and identifying one vertex from each triangle into a single common vertex, called the \emph{center}, such that each triangle shares exactly this center vertex, yielding \(2n + 1\) vertices and \(3n\) edges.

The total vertex irregularity strengths of helm graphs and friendship graphs were first computed by Ahmad et al.~\cite{ahmadTotalVertexIrregularity} and Wijaya et al.~\cite{wijaya2008total}, respectively.

\begin{thm}\label{helm}
Let \(H_n\) be the helm graph with \(2n + 1\) vertices, consisting of a center vertex, an \(n\)-cycle, and \(n\) pendant vertices attached to the cycle vertices. Then
\[
\mathrm{tvs}(H_n) = \left\lceil \frac{n+1}{2} \right\rceil.
\]
\end{thm}

\begin{proof}
First, we establish a lower bound. Suppose \(\mathrm{tvs}(H_n) = k\). Consider the \(n\) pendant vertices of \(H_n\), each of degree 1. The smallest possible weight for a pendant vertex is 2 (when its incident edge and vertex are labeled 1), and the largest is \(2k\) (when labeled \(k\)). Since all \(n\) pendant vertices must have distinct weights, ranging from 2 to \(2k\), we need
\[
n \leq 2k - 1,
\]
so
\[
k \geq \left\lceil \frac{n+1}{2} \right\rceil.
\]
Thus, \(\mathrm{tvs}(H_n) \geq \left\lceil \frac{n+1}{2} \right\rceil\).

To achieve this value, let \(s=\left\lceil \frac{n+1}{2}\right\rceil\). Define a \(\{1,s\}\)-edge labeling on \(H_n\) by labeling every spoke and exactly \(s\) pendant edges with \(1\), and labeling all remaining edges, namely the other \(n-s\) pendant edges and all \(n\) cycle edges, with \(s\).

Then each pendant vertex has vertex sum \(1\) or \(s\), so we obtain \(s\) vertices in \(V_1\) and \(n-s\) vertices in \(V_s\). Every cycle vertex is incident with two cycle edges labeled \(s\) and one spoke labeled \(1\), so its partial sum from non-pendant edges equals \(2s+1\). Consequently, if its pendant edge is labeled \(1\), then its full vertex sum is \(2s+2\), whereas if its pendant edge is labeled \(s\), then its full vertex sum is \(3s+1\). Hence
\[
D_\lambda = (|V_1|, |V_s|, |V_{2s+2}|, |V_{3s+1}|) = (s, n - s, s, n - s),
\]
summing to \(2n\), the number of non-central vertices in \(H_n\). Extend \(\lambda\) to the non-central vertices by
\[
\lambda(\bv{V_1}) = \lambda(\bv{V_{2s+2}}) = \bv{[1,s]}, \quad \lambda(\bv{V_s}) = \lambda(\bv{V_{3s+1}}) = \bv{[2,n+1-s]}.
\]
By the definition of \(s\), all these labels are at most \(s\).

The weights of non-central vertices are
\begin{align*}
\wt(\bv{V_1}) &= \bv{1} + \bv{[1,s]} = \bv{[2,s+1]}, \\
\wt(\bv{V_s}) &= \bv{s} + \bv{[2,n+1-s]} = \bv{[s+2,n+1]}, \\
\wt(\bv{V_{2s+2}}) &= \bv{2s+2} + \bv{[1,s]} = \bv{[2s+3,3s+2]}, \\
\wt(\bv{V_{3s+1}}) &= \bv{3s+1} + \bv{[2,n+1-s]} = \bv{[3s+3,n+2s+2]}.
\end{align*}
By the definition of \(s\), we have \(n+1-s\le s\), so all assigned vertex labels are at most \(s\). Moreover, these intervals are disjoint since \(n+1\le 2s<2s+3\). Since every spoke is labeled \(1\), the center vertex \(x\) has \(w_E(x)=n\). Thus assigning \(\lambda(x)=2\) yields \(\wt(x)=n+2\), which is distinct from all other weights because \(n+1<n+2<2s+3\).
\end{proof}
\begin{thm}\label{friendship}
Let \(F_n\) be the friendship graph with \(2n + 1\) vertices, consisting of a center vertex and \(n\) triangles sharing this center, for \(n \geq 2\). Then
\[
\mathrm{tvs}(F_n) = \left\lceil \frac{2n+2}{3} \right\rceil.
\]
\end{thm}

\begin{proof}
Since \(F_1\) is isomorphic to \(C_3\), with \(\mathrm{tvs}(C_3) = 2\) from Theorem \ref{cycle}, we assume \(n \geq 2\).

First, we establish a lower bound. Suppose \(\mathrm{tvs}(F_n) = k\). Consider the \(2n\) non-central vertices of \(F_n\), each incident to two edges. The smallest possible weight for a non-central vertex is 3  and the largest is \(3k\). Since all \(2n\) non-central vertices must have distinct weights, ranging from 3 to \(3k\), we need
\[
2n \leq 3k - 2,
\]
so
\[
k \geq \left\lceil \frac{2n+2}{3} \right\rceil.
\]
Thus, \(\mathrm{tvs}(F_n) \geq \left\lceil \frac{2n+2}{3} \right\rceil\).

To realize this value, we implement our \(\{1,s\}\)-labeling on \(F_n\), distributing vertex sums evenly across its structure. Let \(s = \left\lceil \frac{2n+2}{3} \right\rceil\) and define a \(\{1,s\}\)-edge labeling on \(F_n\). We classify the edge labelings of the \(n\) triangles as type I, type II, or type III, corresponding to the first, second, and third triangles from left to right in the figure below.  Here the red vertices represent the common center (depicted separately for clarity but forming a single vertex).

\begin{figure}[h]
  \centering 
\includegraphics[width=0.6\textwidth]{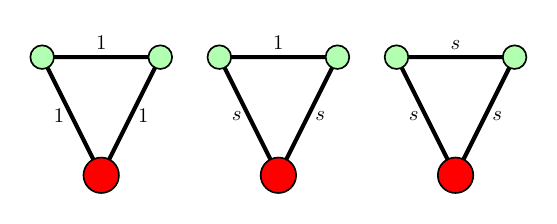} 
  \caption{Type I, II and III edge labelings}
\end{figure}

Let \(V_2\), \(V_{s+1}\), and \(V_{2s}\) be the sets of all non-central vertices of \(F_n\) with vertex sums 2, \(s + 1\), and \(2s\), respectively. Notice that labeling a single triangle using type I, type II, or type III contributes exactly 2 vertices to \(V_2\), \(V_{s+1}\), and \(V_{2s}\), respectively.

If \(s = 2a\) is even, label \(a\) triangles and \(a - 1\) triangles using type I and type II, respectively; then label the rest by using type III. It follows that \(|V_2| = s\), \(|V_{s+1}| = s - 2\), and \(|V_{2s}| = 2n + 2 - 2s\). Then label
\[
\lambda(\bv{V_2}) = \bv{[1,s]}, \quad \lambda(\bv{V_{s+1}}) = \bv{[2,s-1]}, \quad \lambda(\bv{V_{2s}}) = \bv{[1,2n+2-2s]}.
\]
One can easily check that the weights of non-central vertices are all distinct in the range \([3, 2n + 2]\). Lastly, each label is at most \(s\) since
\[
2n + 2 - 2s \leq s \iff 2n + 2 \leq 3s,
\]
which is true by the definition of \(s\).

If \(s\) is odd, then \(s - 1 = 2a\) for some \(a\). Label using type I, type II, and type III accordingly so that
\[
|V_2| = s - 1, \quad |V_{s+1}| = s - 1, \quad \text{and} \quad |V_{2s}| = 2n + 2 - 2s.
\]
Then, similar to the case when \(s\) is even, the labeling
\[
\lambda(\bv{V_2}) = \bv{[1, s - 1]}, \quad \lambda(\bv{V_{s + 1}}) = \bv{[1, s - 1]}, \quad \lambda(\bv{V_{2s}}) = \bv{[1, 2n + 2 - 2s]}
\]
ensures that the weights of all non-central vertices are distinct, and each label is at most \(s\).

Note that at least one triangle is labeled using either type II or type III. This is clear because, otherwise, all non-central vertices would belong to \(V_2\) or \(V_{2s}\). However, this is not possible, since for \(n \geq 2\) we have \(|V_2| \leq s < 2n\).

Consequently, the center vertex is incident to at least two edges labeled \(s\). Since \(F_n\) contains at least two triangles when \(n \geq 2\), assigning the label \(s\) to the center vertex ensures that its weight is at least \(3s + 2\). As the weight of any non-central vertex is at most \(3s\), the weight of the center vertex is distinct from those of the non-central vertices.
\end{proof}
\section{Simple 2-Regular Graphs}
The total vertex irregularity strength of simple 2-regular graphs was previously conjectured by Ahmad et al. \cite{ahmadIrregularTotalLabeling}. They determined the total vertex irregularity strength for the disjoint union of cycles, provided that each cycle has a length of at least 5. However, they were unable to compute it for the general case of simple 2-regular graphs.

\begin{thm}
Let $G$ be a simple $2$-regular graph of order $n$. Then
\[
\tvs(G)=\left \lceil \frac{n+2}{3}\right\rceil.
\]
\end{thm}

\begin{proof}
The graph $G$ is a disjoint union of cycles, $G=\bigcup_{i=1}^m C_i$, where each cycle has length at least $3$. Since the case $m=1$ corresponds to a single cycle (Theorem~\ref{cycle}), we assume $m\ge 2$, which implies $n\ge 6$. Without loss of generality, assume the cycle lengths are ordered such that $n_1 \le n_2 \le \dots \le n_m$, where $n_i = |V(C_i)|$.

Let $s=\lceil (n+2)/3\rceil$. We first consider the case where $s$ is odd. We construct a $\{1,s\}$-edge labeling such that the partition sets $V_2$, $V_{s+1}$, and $V_{2s}$ (sets of vertices with vertex sums $2$, $s+1$, and $2s$, respectively) have sizes
\[
|V_2|=s, \quad |V_{s+1}|=s-1, \quad |V_{2s}|=n+1-2s.
\]
Let $k\ge 0$ be the maximal integer such that $\sum_{i=1}^{k} n_i \le s$, where the empty sum is interpreted as $0$. We write
\[
\sum_{i=1}^{k} n_i = s-a \qquad \text{for some } a \ge 0.
\]
We construct the edge labeling in stages to fix $|V_2|$ exactly at $s$.

\smallskip
\noindent
\textbf{Case 1: $n_{k+1} \ge a+2$.} \\
Label all edges of the first $k$ cycles $C_1, \dots, C_k$ with $1$. On $C_{k+1}$, label a path of $a+1$ consecutive edges with $1$. The remaining edges of $C_{k+1}$ and all edges of subsequent cycles are initially left unlabeled. The cycles $C_1, \dots, C_k$ contribute exactly $\sum_{i=1}^k n_i = s-a$ vertices to $V_2$. The path of $a+1$ edges labeled $1$ in $C_{k+1}$ contributes exactly $(a+1)-1 = a$ internal vertices to $V_2$. Thus, the total count is $|V_2| = (s-a) + a = s$.

\smallskip
\noindent
\textbf{Case 2: $n_{k+1} = a+1$.} \\
If we were to label all $a+1$ edges of $C_{k+1}$ with $1$, the cycle would contribute $a+1$ vertices to $V_2$, yielding a total of $|V_2| = (s-a) + (a+1) = s+1$. Since replacing a label $1$ with $s$ alters the vertex sum of both endpoints, it reduces $|V_2|$ by exactly $2$, making it impossible to achieve $|V_2|=s$ using $C_{k+1}$ alone.

Instead, we label exactly one edge of $C_{k+1}$ with $s$ and the remaining $n_{k+1}-1$ edges with $1$. Then $C_{k+1}$ contributes
\[
n_{k+1}-2=a-1
\]
vertices of vertex sum $2$, so the current total is
\[
(s-a)+(a-1)=s-1.
\]
The two endpoints of the edge labeled $s$ in $C_{k+1}$ now have vertex sum $s+1$. To recover the missing vertex for $V_2$, we use the next cycle, $C_{k+2}$. Note that
\[
\sum_{i=1}^{k+1} n_i = s+1 =\left\lceil\frac{n+2}{3}\right\rceil+1 < \frac{n+2}{3}  + 2< n \quad \text{for } n \ge 6,
\]
ensuring $C_{k+2}$ exists. On $C_{k+2}$, we label exactly two consecutive edges with $1$. This creates exactly one vertex of vertex sum $2$, namely the vertex between the two edges labeled $1$. Therefore
\[
|V_2|=(s-1)+1=s.
\]
The remaining edges of $C_{k+2}$ are left unlabeled for now, which completes Case 2.

\smallskip
\noindent
Thus, in either case, we have fixed $|V_2|=s$. We now determine the labels for the remaining unlabeled edges so as to obtain $|V_{s+1}| = s-1$.

We now complete the labeling of the current cycle (if it is only partially labeled) and label every subsequent fresh cycle using alternating patterns that produce an even number of vertices with vertex sum $s+1$, while creating no new pair of consecutive edges labeled $1$. We use patterns such as $s,1,s,1,\ldots,s,1$ or $s,1,s,1,\ldots,s,1,s$ on a fresh cycle, and patterns such as $s,\mathbf{1},s,1,\ldots,s,1$ or $s,\mathbf{1},s,1,\ldots,s,1,s$ on a cycle that already contains a block of $1$s (denoted by $\mathbf{1}$). Among the \(n-s\) vertices outside \(V_2\), at most one vertex in each remaining cycle has vertex sum \(2s\). Hence the number of vertices with vertex sum \(s+1\) that is produced is at least
\[
n-s-(m-k).
\]
In any case, we will show that this quantity is at least \(s-1\). We now distinguish two subcases. If \(k=0\), then necessarily \(n_1\ge s+1\). It is impossible that \(m\ge 3\), for otherwise
\[
n=\sum_{i=1}^m n_i \ge n_1+n_2+n_3 \ge 3(s+1)>n,
\]
a contradiction. Hence \(m\le 2\). Since we are assuming \(m\ge 2\), it follows that \(m=2\). Therefore \(n=n_1+n_2\ge 2(s+1)\), and after \(|V_2|=s\) has been fixed, the number of vertices with vertex sum \(s+1\) that can be produced is at least
\[
n-s-2\ge s>s-1.
\]
If \(k\ge 1\), then \(m-k\le m-1\), and since every cycle has length at least \(3\), we have \(m\le \left\lfloor n/3 \right\rfloor\). Hence
\[
n-s-(m-k)\ge n-s-\left\lfloor \frac n3 \right\rfloor +1 \ge s-1.
\]
If the resulting number is larger than \(s-1\), we reduce it by replacing local patterns \(s,1,s\) with \(s,s,s\). This operation decreases \(|V_{s+1}|\) by exactly \(2\) and does not affect \(|V_2|\). Indeed, the only vertices of vertex sum \(s+1\) not arising from a local pattern \(s,1,s\) are the two boundary vertices adjacent to the initial block of edges labeled \(1\). Therefore, whenever \(|V_{s+1}|\ge 4\), at least two vertices of \(V_{s+1}\) arise from some local pattern \(s,1,s\), so such a replacement is possible. Since each replacement decreases \(|V_{s+1}|\) by \(2\), the last step before reaching the target \(|V_{s+1}|=s-1\) occurs when \(|V_{s+1}|=s+1\). As \(n\ge 6\), we have \(s=\left\lceil\frac{n+2}{3}\right\rceil\ge 3\), and hence \(s+1\ge 4\), so this final reduction can still be performed. Repeating this operation as needed, we obtain exactly \(|V_{s+1}|=s-1\), while \(|V_2|\) remains unchanged.

Consequently, every vertex not in \(V_2\) or \(V_{s+1}\) has vertex sum \(2s\). Thus we obtain the distribution
\[
|V_2|=s,\qquad |V_{s+1}|=s-1,\qquad |V_{2s}|=n+1-2s.
\]
This is identical to the distribution in Case~2 of Theorem~\ref{cycle}. Therefore, we assign vertex labels to the ordered sets \(\mathbf{V_2}\), \(\mathbf{V_{s+1}}\), and \(\mathbf{V_{2s}}\) according to the same pattern as in that case. Since this assignment depends only on the distribution of these sets and not on the structure of the underlying graph, it guarantees that all vertex weights are pairwise distinct and bounded by $s$.

It remains to consider the case when $s$ is even. The construction used above to force the size of $V_2$ does not rely on the parity of $s$. Hence we may run the same first-stage procedure with target
\[
|V_2|=s-1.
\]
After this stage is fixed, there remain
\[
n-(s-1)
\]
vertices outside $V_2$, distributed among the remaining $(m-k)$ cycle-components under consideration. Using the same alternating completion patterns as in the odd case, at most one vertex per such component fails to contribute to $V_{s+1}$. Therefore
\[
|V_{s+1}|\ge n-(s-1)-(m-k).
\]
Using the inequality proved above,
\[
n-s-(m-k)\ge s-1,
\]
we get
\[
|V_{s+1}|\ge n-(s-1)-(m-k)=(n-s-(m-k))+1\ge (s-1)+1=s.
\]
Hence we can first obtain $|V_{s+1}|\ge s$. The contribution to $V_{s+1}$ from each cycle is even, so the total number $|V_{s+1}|$ obtained after the alternating completion is even. If this number is larger than $s$, we repeatedly apply the local replacement
\[
(s,1,s)\longrightarrow (s,s,s),
\]
which decreases $|V_{s+1}|$ by exactly $2$ and leaves $|V_2|$ unchanged. As in the odd case, whenever a reduction is needed there are enough local patterns $s,1,s$ available outside the initial block of $1$'s to perform it. Since $s$ is even, after finitely many such steps we obtain exactly $|V_{s+1}|=s$.

Consequently, we obtain an edge labeling such that
\[
|V_2|=s-1,\quad |V_{s+1}|=s,\quad \text{and}\quad |V_{2s}|=n+1-2s.
\]
This distribution coincides with that of Case~1 in Theorem~\ref{cycle}. Although the underlying graphs are different, the final vertex-labeling step depends only on the partition \(V(G)=V_2\sqcup V_{s+1}\sqcup V_{2s}\) and the cardinalities of these three classes. Since these coincide with those in Case~1 of Theorem~\ref{cycle}, we assign vertex labels to the ordered sets \(\mathbf{V_2}\), \(\mathbf{V_{s+1}}\), and \(\mathbf{V_{2s}}\) by the same scheme, and the distinct-weight verification applies verbatim. Therefore,
\[
\tvs(G)=\left\lceil\frac{n+2}{3}\right\rceil.
\]

\end{proof}

\section*{Acknowledgment}

The author thanks Institut Teknologi Bandung (ITB) and the Ministry of Higher Education, Science, and Technology (Kementrian Pendidikan Tinggi, Sains dan Teknologi), Republic of Indonesia, for their support in this research.
The author is also grateful to the anonymous reviewer for the careful reading and detailed suggestions, which significantly improved the clarity and presentation of the paper.

\end{document}